\documentclass[english,10pt]{article}
\usepackage{graphicx}
\usepackage{url}
\usepackage{algorithm}
\usepackage{algorithmic}
\usepackage{makeidx}
\usepackage{multicol}
\usepackage[T1]{fontenc}
\usepackage[latin1]{inputenc}
\usepackage{amssymb}
\usepackage[T1]{fontenc}
\usepackage{graphicx}
\usepackage{verbatim}
\makeatletter
\usepackage{graphicx}
\newcommand{\X}{\mathcal{X}}
\newcommand{\Y}{\mathcal{Y}}

\usepackage{amsmath}
\usepackage{makeidx}         
\usepackage{graphicx}        
\usepackage{multicol}        
\usepackage[bottom]{footmisc}

\makeindex             
\begin{document}
\title{Boosting for Functional Data}
\author{Nicole Kr\"amer\thanks{TU Berlin -- Department of Computer Science and Electrical
Engineering,
Franklinstr. 28/29,
10587 Berlin, Germany \texttt{nkraemer@cs.tu-berlin.de}}}

\maketitle
\begin{abstract}
We deal with the task of supervised learning if the data is of
functional type. The crucial point is the choice of the appropriate
fitting method (also called learner).  Boosting is a stepwise
technique that combines learners in such a way that the composite --
boosted -- learner outperforms the single learner. This can be done
by either reweighting the examples or with the help of a gradient
descent technique.  In this paper, we explain how to extend Boosting
methods to problems that involve functional data.
\end{abstract}

{\textbf{Keywords}}: Functional Data Analysis, Boosting
\section{A Short Introduction to Boosting}
\label{secboosting}
The task is the following: We  try to estimate a relationship
\begin{eqnarray}
\label{function1}
 F: \X \rightarrow \Y
\end{eqnarray}
based on a finite set $S=\{(x_1,y_1),\ldots,(x_n,y_n)\} \subset
\mathcal{X}\times \mathcal{Y}$ of observations. A popular strategy
is to fix a class of functions $\mathcal{F}$ and to minimize the
empirical risk
\begin{eqnarray}
\label{emprisk}\frac{1}{n}\sum_{i=1} ^n L( y_i, f(x_i))
\end{eqnarray}
over all elements $f\in \mathcal{F}$. Here
\begin{eqnarray}
\label{loss}
L:\mathcal{Y}\times \mathcal{Y}&\rightarrow &\mathbb{R}
\end{eqnarray}
is a loss function. Sometimes, a regularization term $r(f)$ is added
to (\ref{emprisk}). We call fitting methods like this learners.
Popular examples for multivariate data are trees, support vector
machines or smoothing splines. The choice of the learner is crucial,
as too complex learners lead to overfitting, whilst 'weak' learners
fail to  capture the relevant structure. The term weak learner has
its seeds in the machine learning literature. In classification
problems, a weak learner is a learner that is slightly better than
random guessing. (The exact definition can be found in e.g.
\cite{Meir0301}.) For regression problems, we might think of a
learner that has a high bias compared to its variance, or a learner
that has only a few degrees
of freedom.\\
The basic idea  of Boosting is to proceed stepwise and to combine
weak learners in such a way that the composite -- boosted --
learner
\begin{eqnarray}
\label{combination}
f_M(x)&=& \sum_{m=1} ^M \alpha_m \cdot g_m (x)
\end{eqnarray}
(or $sign(f)$ for classification problems) performs better than the
single weak learners $g_m$.  The single learners are usually called
base learners and $M$ is called the number of Boosting iterations.
The learners $g_m$ and the weights $\alpha_m$ are chosen adaptively
from the data. AdaBoost \cite{Freund9701} -- the first Boosting
algorithm -- is designed for classification problems. It is
presented in algorithm \ref{algoadaboost}. The weak base learner is
repeatedly applied to the weighted training sample $S\,$. Points
which were hard to approximate in step $m$ are given higher weights
in the next iteration step.
\begin{algorithm}
\begin{algorithmic}

\label{algoadaboost}
\STATE{Input: sample $S$, weak learner, $M$, initial weights $D_1(x_i)=1/n$}
\FOR{$m=1,\ldots,M$}
\STATE{Fit a function $g_m$ to the weighted sample $(S,D_m)$
            using the weak learner}
\STATE{Compute the weighted error
            \begin{eqnarray*}
            \epsilon_m=\sum_{i=1} ^n D_m(x_i) I_{\{y_i\not=g_m(x_i)\}}\,.
            \end{eqnarray*}}
\STATE{Set $\alpha_m=\ln \left(\frac{1-\epsilon_m}{\epsilon_m}\right)\,.$}
\STATE{Update the weights:
            \begin{eqnarray*}
            D_{m+1}(x_i)=D_m(x_i) \exp \left(-\alpha_m y_i g_m(x_i)\right)\,.
            \end{eqnarray*}}
\ENDFOR
\RETURN{$h(x)=\text{sign} \left(\sum_{m=1} ^M \alpha_m g_m(x)\right)\,.$}
\end{algorithmic}
\caption{AdaBoost}

\end{algorithm}
For some learners, it  is not possible to compute a weighted loss. Instead, in each step we draw with replacement a sample of
size $n$ from $S$ and use the weights $D_m$ as probabilities.\\
It can be shown \cite{Breiman9801,Breiman9901} that Boosting is a
forward stage-wise fitting method using gradient descent techniques.
More precisely, in each step we fit a weak learner to $x_i$ and the
negative gradient
\begin{eqnarray}
\label{gradient}
u_i&=&-\left. \frac{\partial L(y_i,f)}{\partial
f}\right|_{f=f_m(x_i)}
\end{eqnarray}
of the loss function (\ref{loss}).  The connection between Boosting
and gradient descent methods has lead to a wide range of new algorithms \cite{Friedman0101},
notably  for
regression problems. Note that if we use the quadratic loss
\begin{eqnarray*}
L(y,y')&=& \frac{1}{2} \left(y-y'\right)^2\,,
\end{eqnarray*}
 the negative gradient is simply the vector of
residuals, i.e. we iteratively fit the residuals using a weak
learner. This method  is called $L_2$Boost
\cite{Buehlmann0201,Friedman0101} and is presented in algorithm
\ref{algol2boost}.
\begin{algorithm}
\caption{$L_2$Boost}
\begin{algorithmic}
\label{algol2boost} \STATE{Input: sample $S$, weak learner, $M$}
\STATE{Fit a function $g_1(x)$,   using the weak learner and set
$f_1=g_1$.} \FOR{m=1,\ldots,M} \STATE{Compute the residuals
\begin{eqnarray*}
u_i=& y_i  -  f_m(x_i)\,,i=1,\ldots,n\,.
\end{eqnarray*}}
\STATE{Fit a function $g_{m+1}$ to  $(x_i,u_i)$  by using a weak learner}
\STATE{Update
\begin{eqnarray*}
f_{m+1}=f_m(x)+ g_{m+1}(x)\,.
\end{eqnarray*}}
\ENDFOR
\RETURN{$f_{M}=\sum_{m=1} ^M  g_m (x)$}
\end{algorithmic}
\end{algorithm}
Boosting with the loss function
\begin{eqnarray*}
L(y,y')&=&log\left(1+exp(-yy')\right)\,,
\end{eqnarray*}
is suited for classification problems and called LogitBoost
\cite{Friedman0001}(see algorithm \ref{algologitboost}).
\begin{algorithm}
\caption{LogitBoost}
\begin{algorithmic}
\label{algologitboost} \STATE{Input: sample $S$, weak learner, $M$}
\STATE{Initialize probabilities $p_1(x_i)=1/2$ and set $f_0(x)=0$}
\FOR{m=1,\ldots,M} \STATE{Compute weights  and negative gradients
\begin{eqnarray*}
D_m(x_i)&=& p_m(x_i)\left(1-p_m(x_i)\right)\\
 u_i&=& \frac{y_i-p_{m}(x_i)}{D_m(x_i)},i=1,\ldots,n\,.
\end{eqnarray*}}
\STATE{Fit a regression function $g_{m}$ to  $(x_i,u_i)$  by
weighted least squares} \STATE{Update
\begin{eqnarray*}
f_{m}&=&f_{m-1}(x)+ \frac{1}{2}g_{m}(x)\,,\\
p_{m+1}(x_i)&=& \left(1+exp\left/-2 f_{m}(x_i)\right)\right){-1}\,.
\end{eqnarray*}}
\ENDFOR \RETURN{$f_{M}$}
\end{algorithmic}
\end{algorithm}
The function $f_M$ is an estimate of one-half of the log-odds ratio
\begin{eqnarray*}
\frac{1}{2}\log \left(\frac{P(Y=1|X=x)}{1-P(Y=1|X=x)}\right)\,.
\end{eqnarray*}
As a consequence, this classification algorithm also produces
estimates of the class probabilities $P(Y=1|X=x)$.
 Generic Boosting algorithms for
a general loss function can be found in
\cite{Buehlmann0201,Friedman0101}.\\
How do we obtain the optimal number of Boosting iterations? One possibility
is to use cross validation. Depending on the data, this can lead to high
computational costs. If we use $L_2$Boost, it is possible to compute the
degrees of freedom of the Boosting algorithm
\cite{Buehlmann0201,Buehlmann0501}. As a consequence, we can use model
selection criteria as the Akaike Information Criterion (AIC) or the Bayesian
Information Criterion (BIC).
\section{Functional Data Analysis}
\label{secfda}
The content of this section is condensed from \cite{Ramsay0501}. We speak of
functional data if the variables that we observe are curves. Let us first
consider the case that only the predictor samples $x_i$ are curves, that is
\begin{eqnarray*}
x_i \in \mathcal{X}=\{x:T \rightarrow \mathbb{R}\}\,.
\end{eqnarray*}
Examples for this type of data are time series, temperature curves
or near infra red spectra. We usually assume that the functions
fulfill a regularity condition, and in the rest of the paper, we
consider the Hilbert space $\mathcal{X}=L^2(T)$ of all
square-integrable functions $T\rightarrow \mathbb{R}$.
\subsection{How to Derive Functions from Observations?}
In most applications, we do not measure a curve, but discrete
values of a
curve. An important step in the analysis of functional data is therefore
the transformation of the discretized objects to smooth functions. The
general approach is the following: We represent each example as a linear combination
\begin{eqnarray}
\label{baseexp}
x_i (t)&=& \sum_{l=1} ^{K_x} c_{il} \psi_l (t)
\end{eqnarray}
of a set of base functions $\psi_1,\ldots,\psi_{K_x}$.
The coefficents $c_{il}$ are then estimated by using (penalized) least squares. The most frequently used base functions are Fourier expansions,  B-splines,  wavelets and polynomials. A different possibility is to
derive an orthogonal basis directly from the data. This can be done by using
functional principal component analysis.
\subsection{Inference from Functional Data}
We only consider linear relationships (\ref{function1}), i.e. in the
regression setting ($\mathcal{Y}=\mathbb{R}$), elements $f\in \mathcal{F}=L^2(T)$ are assumed to be linear (up
to an intercept) and
continuous. As $\mathcal{F}$ is a Hilbert space, it follows that any function $f
\in \mathcal{F}$ is of the form
\begin{eqnarray}
\label{linmod}
f(x(t))&=& \beta_0  + \int_{T} \beta(t) x (t) dt \,.
\end{eqnarray}
In the two-class classification setting ($\mathcal{Y}=\{\pm 1\}$), we use  $sign(f)$ instead of $f$. As already mentioned in Sect. \ref{secboosting}, we estimate $f$ or
$\beta$ by minimizing the empirical risk (\ref{emprisk}). Note that this is
an ill-posed problem, as there are (in general) infinitely many
functions $\beta$ that fit the data perfectly. There is obviously a
need for regularization, in order to avoid overfitting. We can solve this problem by using a base expansion of both the predictor
variable $x_i(t)$ as in (\ref{baseexp}) and the function
\begin{eqnarray}
\label{betaexp}
\beta(t)&=& \sum_{l=1} ^{K_{\beta}} b_l \psi_l (t)\,.
\end{eqnarray}
This transforms (\ref{emprisk}) into a parametric problem. If we use the
quadratic loss, this is a matrix problem: We set
\[
\begin{array}{ccccc}
C=(c_{ij})&,&J= \left(\int_T \psi_i(t) \psi_j (t) dt \right) &,& Z=CJ\,.
\end{array}
\]
It follows that (for centered data)
\begin{eqnarray}
\label{bhat}
\hat {\vec{b}} &=& \left(Z^t Z \right)^{-1} Z^t y\,.
\end{eqnarray}
As already mentioned, we have to regularize this problem. There are
two possibilities: We can either constrain the number of base
functions in (\ref{betaexp}). That is, we demand that $K_{\beta}\ll
K_x$. However, we show in Sect. \ref{secfuncboost} that this
strategy can lead to trivial results in the Boosting setting.  The
second possibility is to add a penalty term $r(f)$ to the empirical
risk (\ref{emprisk}). If we consider functional data, it is common
to use a penalty term of the form
\begin{eqnarray*}
r(\beta)&=& \lambda \int_T \left( \beta^{(k)} (t) \right)^2 dt\,.
\end{eqnarray*}
Here $\beta^{(k)}$ is the $k$th derivative of $\beta$ -- provided that this
derivative exists.  The choice of $k$ depends on the data at hand and our expert
knowledge on the problem.\\
Finally, let us briefly mention how to model a linear relationship
(\ref{function1}) if both the predictor and response variable are
functional.  We consider  functions
\begin{eqnarray*}
f:L^2(T)&\rightarrow &L^2(T)\\
f(x(t))&=& \alpha(t) + \int_T \beta(s,t) x (s) ds \,.
\end{eqnarray*}
We estimate $\beta$ by expanding $y_i, x_i, \alpha$ in terms of a basis and
by representing $\beta$ by
\begin{eqnarray*}
\beta(s,t)&=& \sum_{k=1} ^{K_1} \sum_{l=1}^{K_2} b_{kl} \psi_k(s) \psi_l(t)\,.
\end{eqnarray*}
The optimal coefficients $b_{kl}$ are determined using the loss function
\begin{eqnarray*}
L(y,y')&=& \int_T \left(y(t)-y'(t)\right)^2 dt\,.
\end{eqnarray*}
Again, we have to regularize in order to obtain smooth estimates that do not
overfit.
\section{Functional Boosting}
\label{secfuncboost}
In order to apply a Boosting technique to functional data, we have to extend
the notion 'weak learner'. In the classification setting, we can adopt the loose definition
from Sect. \ref{secboosting}. A weak learner is a learner that is
slightly better than random. What are examples of weak learners?
Note that it is possible to apply most of the multivariate data
analysis tools to functional data. We use a finite-dimensional
approximation as in (\ref{baseexp}) and simply apply any appropriate
algorithm. In this way, it is possible to use stumps (that is,
classification trees with one node) or neural networks as base
learners.\\
In the regression setting, we propose the following definition:  A
weak learner is a learner that has only a few degrees of freedom.
Examples include the two regularized least squares algorithms
presented in Sect. \ref{secfda} --  restriction of the number of base
functions in (\ref{betaexp}) or addition of a penalty term to
(\ref{emprisk}). Note however that the first method leads to trivial
results if we use $L_2$Boost. The learner is simply the projection
of $y$ onto the space that is spanned by the columns of $Z$ (recall
(\ref{bhat})). Consequently, the $y$-residuals are orthogonal on $Z$
and after one step, the Boosting solution does not change anymore.
Another example of a weak learner is the following
\cite{Buehlmann0501}: In each Boosting step, we only select one base
function using $x_i$ and the residuals $u_i$. To select this base
function, we estimate the regression coefficients $b_j$ of
\begin{eqnarray}
\label{sparse}
u_i &\sim& b_j \int_T x_i(t) \psi_j(t) dt\,,\,j=1,\ldots,K_{\beta}\,.
\end{eqnarray}
We choose the base function that
minimizes the empirical risk (\ref{emprisk}). For centered data, this equals
\begin{eqnarray*}
m^*&=&\text{arg}\min_{j} \sum_{i=1} ^n L\left(u_i,\hat b_j \int_T
\psi_j(t)x_i(t)dt \right)\\
\hat b_j&=& \text{Least Squares estimate of }b_j\text{ in (\ref{sparse})}
\end{eqnarray*}
Boosting for multivariate data with this kind of weak learner has been studied in
e.g. \cite{Buehlmann0501}. \\
If the response variable is functional, we can adopt the same
definition of weak learner as in the regression setting: A weak
learner is a learner that uses only a few degrees of freedom.
\section{Example: Speech Recognition}
This example is taken from \cite{Biau0501}. The data consists of
$48$ recordings  of the word 'Yes' and $52$ recordings of the word
'No'. One recording is represented by a discretized time series of
length $8192$. The data can be downloaded from
\url{http://www.math.univ-montp2.fr/~biau/bbwdata.tgz}. All
calculations
are performed using \texttt{R} \cite{R}. \\
The task is to find a classification rule that assigns the correct
word to each time series. We apply the LogitBoost algorithm to this
data set. First, we represent the time series in terms of a Fourier
basis expansion of dimension $K_x=100$. We opted to include a lot of basis
functions, as experiments indicate that the results of
LogitBoost are insensitive to the addition of possibly irrelevant  basis functions. The weak
learner is a classification tree with two final nodes. The
misclassification rate was estimated using 10fold cross-validation
(cv).
\begin{figure}[h]
{\par\centering\resizebox*{8cm}{5.5cm}{\rotatebox{270}{{\includegraphics{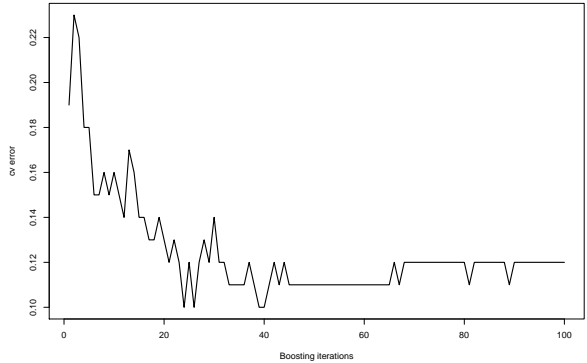}}}}\par}
\caption{Cross validated error for the speech recognition problem.
The optimal number of Boosting iterations is $M_{opt}=24$. }
\label{cverror}
\end{figure}
Figure \ref{cverror} shows the cross-validated error as a function
of the number of Boosting iterations. The minimal cv error over all
Boosting iterations is 0.1, obtained after 24 Boosting iterations.
This is the same error rate that is reported in \cite{Biau0501}.
There, a functional $k$-nearest-neighbor-algorithm is applied to the data.
Finally, we remark that the cv error curve stays rather flat after
the minimum is attained. This  seems to be a feature of all Boosting
methods. As a consequence, the selection of the optimal number of
Boosting iterations can be done quite easily.
\section{Conclusion}
The extension of Boosting methods to functional data is
straightforward. After choosing a base algorithm (which we called a weak learner), we iteratively  fit the data by either applying
this algorithm to reweighted samples or by using a gradient descent
technique. In many applications, we use a finite-dimensional
expansion of the functional examples in terms of base functions.
This finite-dimensional representation can then be plugged into
existing algorithms as AdaBoost, LogitBoost or
$L_2$Boost. \\
We focused on linear learning problems in Sect. \ref{secfda} for the
sake of simplicity and briefness, but it should be noted that
Boosting methods can also be applied to solve nonlinear functional
data problems.

\end{document}